\documentclass[10pt,twoside]{article}
\usepackage{makeidx}
\usepackage{amsfonts}
\usepackage{amsmath}
\usepackage{geometry}
\usepackage{euler}
\usepackage{times}

\newtheorem{theorem}{Theorem}

\newtheorem{condition}[theorem]{Condition}

\newtheorem{definition}[theorem]{Definition}
\newtheorem{example}[theorem]{Example}

\newtheorem{lemma}[theorem]{Lemma}

\newtheorem{proposition}[theorem]{Proposition}

\input{tcilatex}
\begin{document}

\title{A Proof of a Non-Commutative Central Limit Theorem by the Lindeberg
Method}
\author{Vladislav Kargin \thanks{%
Courant Institute of Mathematical Sciences; 109-20 71st Road, Apt. 4A,
Forest Hills NY 11375; kargin@cims.nyu.edu}}
\date{}
\maketitle

\begin{abstract}
A Central Limit Theorem for non-commutative random variables is proved using
the Lindeberg method. The theorem is a generalization of the Central Limit
Theorem for free random variables proved by Voiculescu. The Central Limit
Theorem in this paper relies on an assumption which is weaker than freeness.
\end{abstract}

\section{Introduction}

One of the most important results in free probability theory is the Central
Limit Theorem (CLT) for free random variables (\cite{voiculescu83}). It was
proved almost simultaneously with the invention of free probability theory.
Later conditions of the theorem were relaxed (\cite{voiculescu98}).
Moreover, a far-reaching generalization was achieved in \cite%
{bercovici_pata99}, which studied domains of attraction of probability laws
with respect to free additive convolutions. See also \cite%
{chistyakov_gotze06}.

Freeness is a very strong condition imposed on operators and it is of
interest to find out whether the Central Limit Theorem continues to hold if
this condition is somewhat relaxed. This problem calls for a different proof
of the non-commutative CLT which does not depend on $R$-transforms or on the
vanishing of mixed free cumulants, because both of these techniques are
closely connected with the concept of freeness.

In this paper we give a proof of free CLT that avoids using either $R$%
-transforms or free cumulants. This allows us to develop a generalization of
the free CLT to random variables that are not necessarily free but that
satisfy a weaker assumption. An example shows that this assumption is
strictly weaker than the assumption of freeness.

The proof that we use is a modification of the Lindeberg proof of the
classical CLT (\cite{lindeberg22}). The main difference is that we use
polynomials instead of arbitrary functions from $C_{c}^{3}\left( \mathbb{R}%
\right) ,$ and that more ingenuity is required to estimate the residual
terms in the Taylor expansion formula.

The closest result to the result in this paper is Theorem 2.1 in (\cite%
{voiculescu91}), where the Central Limit Theorem is proved under the
conditions on summands that are weaker than the requirement of freeness. The
conditions that we use are somewhat different than those in Voiculescu's
paper. In addition, we give an explicit example of variables that are not
free but that satisfy conditions of the theorem.

The rest of the paper is organized as follows. Section 2 provides background
material and formulates the main result. Section 3 shows by an example that
a condition in the main result is strictly weaker than the condition of
freeness. Section 4 contains the proof of the main result. And Section 5
concludes.

\section{Background and Main Theorem}

Before proceeding further, let us establish the background. A \emph{%
non-commutative random space} $\left( \mathcal{A},E\right) $ is a pair of an
operator algebra $\mathcal{A}$ and a linear functional $E$ on $\mathcal{A}$.
It is assumed that $\mathcal{A}$ is closed relative to taking the adjoints
and contains a unit, and that $E$ is

1) positive, i.e., $E\left( X^{\ast }X\right) \geq 0$ for every $X\in 
\mathcal{A}$,

2) finite, i.e., $E(I)=1$ where $I$ denotes the unit operator, and

3) tracial, i.e., $E\left( X_{1}X_{2}\right) =E\left( X_{2}X_{1}\right) $
for every $X_{1}$ and $X_{2}\in \mathcal{A}$.

This linear functional is called \emph{expectation}. Elements of $\mathcal{A}
$ are called \emph{random variables}.

Let $X$ be a self-adjoint random variable (i.e., a self-adjoint operator
from algebra $\mathcal{A}$). We can write $X$ as an integral over a
resolution of identity:%
\begin{equation*}
X=\int_{-\infty }^{\infty }\lambda dP_{X}\left( \lambda \right) ,
\end{equation*}%
where $P_{X}\left( \lambda \right) $ is an increasing family of commuting
projectors. Then we can define the \emph{spectral probability measure} 
\index{spectral probability measure|textbf}of interval $\left( a,b\right] $
as follows:%
\begin{equation*}
\mu _{X}\left\{ \left( a,b\right] \right\} =E\left[ P_{X}\left( b\right)
-P_{X}\left( a\right) \right] .
\end{equation*}%
We can extend this measure to all measurable subsets in the usual way. We
will call $\mu _{X}$ the \emph{spectral probability measure} of random
variable $X,$ or simply its \emph{spectral measure}.

We can calculate the expectation of any summable function of a self-adjoint
variable $X$ by using its spectral measure:%
\begin{equation*}
Ef\left( X\right) =\int_{-\infty }^{\infty }f\left( \lambda \right) d\mu
_{X}\left( \lambda \right) .
\end{equation*}%
In particular, the\emph{\ moments} of the probability measure $\mu _{X}$
equal the expectation values of the powers of $X$:%
\begin{equation*}
\int_{-\infty }^{\infty }\lambda ^{k}d\mu _{X}\left( \lambda \right)
=E\left( X^{k}\right) .
\end{equation*}

Let us now recall the definition of freeness. Consider sub-algebras $%
\mathcal{A}_{1,}...,\mathcal{A}_{n}$. Let $a_{i}$ denote elements of these
sub-algebras and let $k\left( i\right) $ be a function that maps the index
of an element to the index of the corresponding algebra: $a_{i}\in \mathcal{A%
}_{k\left( i\right) }.$

\begin{definition}
The algebras $\mathcal{A}_{1,}...,\mathcal{A}_{n}$ (and their elements) are 
\emph{free} if $E\left( a_{1}...a_{m}\right) =0$ whenever the following two
conditions hold:\newline
(a) $E\left( a_{i}\right) =0$ for every $i$, and \newline
(b) $k(i)\neq k\left( i+1\right) $ for every $i<m.$\newline
The variables $X_{1},...,X_{n}$ are called free if the algebras $\mathcal{A}%
_{i}$ generated by $\left\{ I,X_{i},X_{i}^{\ast }\right\} ,$ respectively,
are free.
\end{definition}

An important property of freeness is that we can compute the moments of the
products of the free random variables.

\begin{proposition}
\label{joint_moment_reduction} Suppose $X_{1,}...,X_{n}$ are free. Then 
\begin{equation}
E\left( X_{1}...X_{n}\right) =\sum_{r=1}^{n}\sum_{1\leq k_{1}<...<k_{r}\leq
n}\left( -1\right) ^{r-1}E\left( X_{k_{1}}\right) ...E\left(
X_{k_{r}}\right) E\left( X_{1}...%
\widehat{X}_{k_{1}}...\widehat{X}_{k_{r}}...X_{n}\right) ,
\label{free_relation_property1}
\end{equation}%
where $\symbol{94}$ denotes terms that are omitted.\newline
\end{proposition}

This property is easy to prove by induction. However, we will not need all
the power of this property. Below we formulate the conditions that we need
to impose on the random variables to prove the CLT. These conditions are
consequences of freeness but are likely to be weaker.

We will say that a sequence of \emph{zero-mean} random variables $%
X_{1},...,X_{n},...$ satisfies \emph{Condition A} if:

\begin{enumerate}
\item For every $k,$ $E\left( X_{k}X_{i_{1}}...\text{ }X_{i_{r}}\right) =0$
provided that $i_{s}\neq k$ for $s=1,...,r.$

\item For every $k\geq 2,$ $E\left( X_{k}^{2}X_{i_{1}}...\text{ }%
X_{i_{r}}\right) =E\left( X_{k}^{2}\right) E\left( X_{i_{1}}...\text{ }%
X_{i_{r}}\right) $ provided that $i_{s}<k\,\ $for $s=1,...,r.$

\item For every $k\geq 2,$ 
\begin{equation*}
E\left( X_{k}X_{i_{1}}...\text{ }X_{i_{p}}X_{k}X_{i_{p+1}}...\text{ }%
X_{i_{r}}\right) =E\left( X_{k}^{2}\right) E\left( X_{i_{1}}...\text{ }%
X_{i_{p}}\right) E\left( X_{i_{p+1}}...\text{ }X_{i_{r}}\right)
\end{equation*}%
provided that $i_{s}<k\,\ $for $s=1,...,r.$
\end{enumerate}

Intuitively, if we know how to calculate every moment of the sequence $%
X_{1},...,X_{k-1},$ then using Condition A we can also calculate the
expectation of any product of random variables $X_{1},$ $...,$ $X_{k}$ that
involves no more than two occurrences of variable $X_{k}$. Part 1 of
Condition A is stronger than is needed for this calculation, since it
involves variables with indices higher than $k.$ However, we will need this
additional strength in the proof of Lemma \ref{moment_bound0} below, which
is essential for the proof of the main result.

\begin{proposition}
Every sequence of free random variables $X_{1},$ $...,$ $X_{n}$, $...$
satisfies Condition A.
\end{proposition}

This proposition can be checked by direct calculation using Proposition \ref%
{joint_moment_reduction}.

We will also need the following fact.

\begin{proposition}
\label{proposition_A1}Let $X_{1},...,X_{l}$ be zero-mean variables that
satisfy Condition A(1), and let $Y_{l+1},...,Y_{n}$ be zero-mean variables
which are free from each other and from the algebra generated by variables $%
X_{1},...,X_{l}.$ Then the sequence $X_{1},...,X_{l},Y_{l+1},...,Y_{n}$
satisfies Condition A(1).
\end{proposition}

\textbf{Proof:} Consider the moment $E\left(
X_{k}A_{i_{1}}...A_{i_{s}}\right) ,$ where $A_{i_{t}}$ is either one of $%
Y_{j}$ or one of $X_{i}$ but it can equal $X_{k}.$ Then we can use the fact
that $Y_{j}$ are free and write 
\begin{equation*}
E\left( X_{k}A_{i_{1}}...A_{i_{s}}\right) =\sum_{\alpha }c_{\alpha }E\left(
X_{k}X_{i_{1}\left( a\right) }...X_{i_{r}\left( \alpha \right) }\right) ,
\end{equation*}%
where none of $X_{i_{t}\left( \alpha \right) }$ equals $X_{k}.$ Then, using
the assumption that $X_{i}$ satisfy Condition A(1), we conclude that $%
E\left( X_{k}A_{i_{1}}...A_{i_{s}}\right) =0.$ Also, $E\left(
Y_{k}A_{i_{1}}...A_{i_{s}}\right) =E\left( Y_{k}\right) E\left(
A_{i_{1}}...A_{i_{s}}\right) =0,$ provided that none of $A_{i_{t}}$ equals $%
Y_{k}.$ In sum, the sequence $X_{1},...,X_{l},Y_{l+1},...,Y_{n}$ satisfies
Condition A(1). QED.

While the freeness of random variables $X_{i}$ is the same concept as the
freeness of the algebras that they generate, Condition A deals only with
variables $X_{i}$, and not with the algebras that they generate. For
example, it is conceivable that a sequence $\left\{ X_{i}\right\} $
satisfies condition $A$ but $\left\{ X_{i}^{2}-E\left( X_{i}^{2}\right)
\right\} $ does not. In particular, this implies that Condition A requires
checking a much smaller set of moment conditions than freeness. Below we
will present an example of random variables which are not free but which
satisfy Condition A.

Recall that the standard semicircle law $\mu _{SC}$ is the probability
distribution on $\mathbb{R}$ with the density $\pi ^{-1}\sqrt{4-x^{2}}$ if $%
x\in \lbrack -2;2],$ and $0$ otherwise. We are going to prove the following
Theorem.

\begin{theorem}
\label{CLT}Suppose that \newline
(i) $\left\{ \xi _{i}\right\} $ is a sequence of self-adjoint random
variables that satisfies Condition $A;$\newline
(ii) every $\xi _{i}$ has asbsolute moments of all orders, which are
uniformly bounded, i.e., $E\left\vert \xi _{i}\right\vert ^{k}\leq \mu _{k}$
for all $i;$\newline
(iii) $E\xi _{i}=0,$ $E\xi _{i}^{2}=\sigma _{i}^{2}$; \newline
(iv) $\left( \sigma _{1}^{2}+...+\sigma _{N}^{2}\right) /N\rightarrow s$ as $%
N\rightarrow \infty $. \newline
Then the spectral measure of $S_{N}=\left( \xi _{1}+...+\xi _{N}\right) /%
\sqrt{\sigma _{1}^{2}+...+\sigma _{N}^{2}}$ converges in distribution to the
semicircle law $\mu _{SC}.$
\end{theorem}

The contribution of this theorem is twofold. First, it shows that the
semicircle central limit holds for a certain class of non-free variables.
Second, it gives a proof of the free CLT which is different from the usual
proof through $R$-transforms. However, it is not stronger than a version of
the free CLT which is formulated in Section 2.5 in \cite{voiculescu98}.

\section{Example}

Let us present an example that suggest that Condition A is strictly weaker
than the freeness condition.

Let $F$ be the free group with a countable number of generators $f_{k}$.
Consider the set of relations $R=\left\{
f_{k}f_{k-1}f_{k}f_{k-1}f_{k}f_{k-1}=e\right\} ,$ where $k\geq 2,$ and
define $G=F/\mathcal{R},$ that is, $G$ is the group with generators $f_{k}$
and relations generated by relations in $R.$

Here are some consequences of these relationships: \newline
1) $f_{k-1}f_{k}f_{k-1}f_{k}f_{k-1}f_{k}=e.$ \newline
(Indeed, $e=f_{k}^{-1}\left( f_{k}f_{k-1}f_{k}f_{k-1}f_{k}f_{k-1}\right)
f_{k}=f_{k-1}f_{k}f_{k-1}f_{k}f_{k-1}f_{k}.$)\newline
2) $f_{k-1}^{-1}f_{k}^{-1}f_{k-1}^{-1}f_{k}^{-1}f_{k-1}^{-1}f_{k}^{-1}=e$
and $f_{k}^{-1}f_{k-1}^{-1}f_{k}^{-1}f_{k-1}^{-1}f_{k}^{-1}f_{k-1}^{-1}=e.$

\bigskip

We are interested in the structure of the group $G$. For this purpose we
will study the structure of $\mathcal{R}$, which is a subgroup of $F$
generated by elements of $R$ and their conjugates. We will represent
elements of $F$ by \emph{words}, that is, by sequences of generators. We
will say that a word is\emph{\ reduced} if does not have a subsequence of
the form $f_{k}f_{k}^{-1}$ or $f_{k}^{-1}f_{k}.$ It is \emph{cyclically
reduced} if it does not have the form of $f_{k}...f_{k}^{-1}$ or $%
f_{k}^{-1}...f_{k}.$ We will call a number of elements in a reduced word $w$
its length and denote it as $\left\vert w\right\vert .$ A set of relations $%
R $ is \emph{symmetrized} if for every word $r\in $ $R,$ the set $R$ also
contains its inverse $r^{-1}$ and all cyclically reduced conjugates of both $%
r$ and $r^{-1}.$

For our particular example, a symmetrized set of relations is given by the
following list:

\begin{equation*}
R=\left\{ 
\begin{array}{c}
f_{k}f_{k-1}f_{k}f_{k-1}f_{k}f_{k-1},\text{ \ }%
f_{k-1}f_{k}f_{k-1}f_{k}f_{k-1}f_{k},\text{ } \\ 
\text{\ }f_{k-1}^{-1}f_{k}^{-1}f_{k-1}^{-1}f_{k}^{-1}f_{k-1}^{-1}f_{k}^{-1},%
\text{ \ }f_{k}^{-1}f_{k-1}^{-1}f_{k}^{-1}f_{k-1}^{-1}f_{k}^{-1}f_{k-1}^{-1}%
\end{array}%
\right\} ,
\end{equation*}
where $k$ are all integers $\geq 2.$

A word $b$ is called a \emph{piece} (relative to a symmetrized set $R$) if
there exist two elements of $R,$ $r_{1}$ and $r_{2},$ such that $%
r_{1}=bc_{1} $ and $r_{2}=bc_{2}.$ In our case, each $f_{k}$ and $f_{k}^{-1}$
with index $k\geq 2$ is a piece because $f_{k}$ is the initial part of
relations $f_{k}f_{k-1}f_{k}f_{k-1}f_{k}f_{k-1}$ and $%
f_{k}f_{k+1}f_{k}f_{k+1}f_{k}f_{k+1},$ and $f_{k}^{-1}$ is the initial part
of relations $%
f_{k}^{-1}f_{k-1}^{-1}f_{k}^{-1}f_{k-1}^{-1}f_{k}^{-1}f_{k-1}^{-1}$ and $%
f_{k}^{-1}f_{k+1}^{-1}f_{k}^{-1}f_{k+1}^{-1}f_{k}^{-1}f_{k+1}^{-1}.$ There
is no other pieces.

Now we introduce the condition of \emph{small cancellation} for a
symmetrized set $R$:

\begin{condition}[$C^{\prime }\left( \protect\lambda \right) $]
If $r\in R$ and $r=bc$ where $b$ is a piece, then $\left\vert b\right\vert
<\lambda \left\vert r\right\vert .$
\end{condition}

Essentially, the condition says that if two relations are multiplied
together, then a possible cancellation must be relatively small. Note that
if $R$ satisfies $C^{\prime }\left( \lambda \right) $ then it satisfies $%
C^{\prime }\left( \mu \right) $ for all $\mu \geq \lambda .$

In our example $R$ satisfies $C^{\prime }\left( 1/5\right) .$

Another important condition is the triangle condition.

\begin{condition}[$T$]
Let $r_{1},$ $r_{2},$ and $r_{3}$ be three arbitrary elements of $R$ such
that $r_{2}\neq r_{1}^{-1}$ and $r_{3}\neq r_{2}^{-1}$ \ Then at least one
of the products $r_{1}r_{2}$, $r_{2}r_{3},$ or $r_{3}r_{1},$ is reduced
without cancellation.
\end{condition}

In our example, Condition $\left( T\right) $ is satisfied.

If $s$ is a word in $F,$ then $s>\lambda R$ means that there exists a word $%
r\in R$ such that $r=st$ and $\left\vert s\right\vert >\lambda \left\vert
r\right\vert .$ An important result from small cancellation theory that we
will use later is the following theorem:

\begin{theorem}[Greendlinger's Lemma]
Let $R$ satisfy $C^{\prime }\left( 1/4\right) $ and $T.$ Let $w$ be a
non-trivial, cyclically reduced word with $w\in \mathcal{R}$. Then either%
\newline
(1) $w\in R,$ \newline
or some cyclycally reduced conjugate $w^{\ast }$ of $w$ contains one of the
following:\newline
(2) two disjoint subwords, each $>\frac{3}{4}R,$ or \newline
(4) four disjoint subwords, each $>\frac{1}{2}R.$
\end{theorem}

This theorem is Theorem 4.6 on p. 251 in \cite{lyndon_schupp77}.

Since in our example $R$ satisfies both $C^{\prime }\left( 1/4\right) $ and $%
T,$ we can infer that in our case the conclusion of the theorem must hold.
For example, (2) means that we can find two disjoint subwords of $w,$ $s_{1}$
and $s_{2},$ and two elements of $R$, $r_{1}$ and $r_{2},$ such that $%
r_{i}=s_{i}t_{i}$ and $\left\vert s_{i}\right\vert >\left( 3/4\right)
\left\vert r_{i}\right\vert =9/2.$ In particular, we can conclude that in
this case $\left\vert w\right\vert \geq 10.$ Similarly, in case (4), $%
\left\vert w\right\vert \geq 16.$ One immediate application is that $G$ does
not collapse into the trivial group. Indeed, $f_{i}$ are not zero.

Let $L^{2}\left( G\right) $ be the functions of $G$ that are square-summable
with respect to the counting measure. $G$ acts on $L^{2}\left( G\right) $ by
left translations: 
\begin{equation*}
\left( L_{g}x\right) \left( h\right) =x\left( gh\right) .
\end{equation*}%
Let $\mathcal{A}$ be the group algebra of $G.$ The action of $G$ on $%
L^{2}\left( G\right) $ can be extended to the action of $\mathcal{A}$ on $%
L^{2}\left( G\right) .$ Define the expectation on this group algebra \ by
the following rule:%
\begin{equation*}
E\left( h\right) =\left\langle \delta _{e},L_{h}\delta _{e}\right\rangle ,
\end{equation*}%
where $\left\langle \cdot ,\cdot \right\rangle $ denotes the scalar product
in $L^{2}\left( G\right) .$ Alternatively, the expectation can be written as
follows:%
\begin{equation*}
E\left( h\right) =a_{e},
\end{equation*}%
where $h=\sum_{g\in G}a_{g}g$ is a representation of a group algebra element 
$h$ as a linear combination of elements $g\in G.$ The expectation is clearly
positive and finite by definition. It is also tracial because $g_{1}g_{2}=e$
if and only if $g_{2}g_{1}=e.$

If $L_{h}=\sum_{g\in G}a_{g}L_{g}$ is a linear operator corresponding to the
element of group algebra $h=\sum_{g\in G}a_{g}g$, then its adjoint is $%
\left( L_{h}\right) ^{\ast }=\sum_{g\in G}\overline{a_{g}}L_{g^{-1}},$ which
corresponds to the element $h^{\ast }=\sum_{g\in G}\overline{a}_{g}g^{-1}.$

Consider elements $X_{i}=f_{i}+f_{i}^{-1}.$ They are self-adjoint and $%
E\left( X_{i}\right) =0.$ Also we can compute $E\left( X_{i}^{2}\right) =2.$
Indeed it is enough to note that $f_{i}^{2}\neq e$ and $f_{i}^{-2}\neq e,$
and this holds because insertion or deletion of an element from $R$ changes
the degree of $f_{i}$ by a multiple of $3.$ Therefore, every word equal to
zero must have the degree of every $f_{i}$ equal to 0 modulo 3.

\begin{proposition}
The sequence of variables $\left\{ X_{i}\right\} $ is not free but satisfies
Condition A.
\end{proposition}

\textbf{Proof:} The variables $X_{k}$ are not free. Consider $%
X_{2}X_{1}X_{2}X_{1}X_{2}X_{1}.$ Its expectation is $2,$ because $%
f_{2}f_{1}f_{2}f_{1}f_{2}f_{1}=e$ and $%
f_{2}^{-1}f_{1}^{-1}f_{2}^{-1}f_{1}^{-1}f_{2}^{-1}f_{1}^{-1}=e,$ and all
other terms in the expansion of $X_{2}X_{1}X_{2}X_{1}X_{2}X_{1}$ are
different from $e.$ Indeed, the only terms that are not of the form above
but still have the degree of all $f_{i}$ equal to zero modulo 3 are $%
f_{2}f_{1}^{-1}f_{2}f_{1}^{-1}f_{2}f_{1}^{-1}$ and $%
f_{2}^{-1}f_{1}f_{2}^{-1}f_{1}f_{2}^{-1}f_{1},$ but they do not equal zero
by application of Greendlinger's lemma. Therefore, $E\left(
X_{2}X_{1}X_{2}X_{1}X_{2}X_{1}\right) =2.$ This contradicts the definition
of freeness of variables $X_{2}$ and $X_{1}.$

Let us check Condition A. For A(1), we have to prove that $E\left(
X_{k}X_{i_{1}}...X_{i_{n}}\right) =0,$ where $k\neq i_{s}$ and $i_{s}\neq
i_{s+1}$ for every $s.$ Consider $E\left( f_{k}f_{i_{1}}...f_{i_{n}}\right)
, $ where $k\neq i_{s}$ and $i_{s}\neq i_{s+1}$ for every $s.$ Note $%
f_{k}f_{i_{1}}...f_{i_{n}}\neq e,$ as can be seen from the fact that the
degree of $f_{k}$ does not equal zero modulo 3. Therefore $E\left(
f_{k}f_{i_{1}}...f_{i_{n}}\right) =0.$ A similar argument works for $E\left(
f_{k}^{-1}f_{i_{1}}...f_{i_{n}}\right) =0$ and more generally for the
expectation of every element of the form $f_{k}^{\varepsilon
}f_{i_{1}}^{n_{1}}...f_{i_{n}}^{n_{2}},$ where $\varepsilon =\pm 1$ and $%
n_{s}$ are integer.

Similarly, we can prove that $E\left( f_{k}^{\pm
2}f_{i_{1}}^{n_{1}}...f_{i_{n}}^{n_{2}}\right) =0$ and this suffices to
prove A(2).

For A(3) we have to consider elements of the form $f_{k}^{\varepsilon
_{1}}f_{i_{1}}...f_{i_{p}}f_{k}^{\varepsilon _{2}}f_{i_{p+1}}...f_{i_{q}}.$
Assume that neither $f_{i_{1}}...f_{i_{p}}$ nor $f_{i_{p+1}}...f_{i_{q}}$
can be reduced to $e.$ Otherwise we can use property A2. Then the claim is
that $E\left( f_{k}^{\varepsilon
_{1}}f_{i_{1}}...f_{i_{p}}f_{k}^{\varepsilon
_{2}}f_{i_{p+1}}...f_{i_{q}}\right) =0.$ This is clear when $\varepsilon
_{1} $ and $\varepsilon _{2}$ have the same sign since in this case the
degree of $f_{k}$ does not equal $0$ modulo $3.$ A more difficult case is
when $\varepsilon _{1}=1$ and $\varepsilon _{2}=-1.$ (The case with opposite
signs is similar.) However, in this case we can conclude that $%
f_{k}f_{i_{1}}...f_{i_{p}}f_{k}^{-1}f_{i_{p+1}}...f_{i_{q}}\neq e$ by an
application of Greendlinger's lemma. Indeed, the only subwords that this
word can contain and which would also be subwords of an element of R, are
subwords of length 1 and 2. But these subwords fail to satisfy the
requirement of either (2) or (4) in Greendlinger's lemma. Therefore, we can
conclude that $f_{k}f_{i_{1}}...f_{i_{p}}f_{k}^{-1}f_{i_{p+1}}...f_{i_{q}}%
\neq e,$ and therefore A(3) is also satisfied. Thus Condition A is satisfied
by random variables $X_{1},...,X_{k},...$ in algebra $\mathcal{A},$ although
these variables are not free. QED.

\section{Proof of the Main Result}

\textbf{Outline of Proof}: Our proof of the free CLT proceeds along the
familiar lines of the Lindeberg method. We take a family of functions, $%
\left\{ f\right\} ,$ and compare $Ef\left( S_{N}\right) $ with $Ef\left( 
\widetilde{S}_{N}\right) ,$ where $S_{N}=X_{1}+...+X_{N}$ and $\widetilde{S}%
_{N}=Y_{1}+...+Y_{N},$ and $Y_{i}$ are free semicircle variables chosen in
such a way that $\mathrm{Var}\left( S_{N}\right) =\mathrm{Var}\left( 
\widetilde{S}_{N}\right) .$ To estimate $\left\vert Ef\left( S_{N}\right)
-Ef\left( \widetilde{S}_{N}\right) \right\vert $, we substitute the elements
in $S_{N}$ with free semicircle variables, one by one, and estimate the
corresponding change in the expected value of $f\left( S_{N}\right) $. After
that, we show that the total change, as all elements in the sum are
substituted with semicircle random variables, is asymptotically small as $%
N\rightarrow \infty .$ Finally, the tightness of the selected family of
functions allows us to conclude that the distribution of $S_{N}$ must
converge to the semicircle law as $N\rightarrow \infty .$

The usual choice of functions $f$ in the classical case are functions from $%
C_{c}^{3}\left( \mathbb{R}\right) ,$ that is, functions with a continuous
third derivative and compact support. In the non-commutative setting this
family of functions is not appropriate because the usual Taylor series
formula is difficult to apply. Intuitively, it is difficult to develop $%
f\left( X+h\right) $ in a power series of $h$ if variables $X$ and $h$ do
not commute. Since the Taylor formula is crucial for estimating the change
in $Ef\left( S_{N}\right) $, we will still use it but we will restrict the
family of functions to polynomials.

To show that the family of polynomials is sufficiently rich for our
purposes, we use the following Proposition:

\begin{proposition}
\label{weak_convergence}Suppose there is a unique distribution function F
with the moments $\left\{ m^{\left( r\right) },r\geq 1\right\} .$ Suppose
that $\left\{ F_{N}\right\} $ is a sequence of distribution functions, each
of which has all its moments finite:%
\begin{equation*}
m_{N}^{\left( r\right) }=\int_{-\infty }^{\infty }x^{r}dF_{N}.
\end{equation*}%
Finally, suppose that for every $r\geq 1:$%
\begin{equation*}
\lim_{n\rightarrow \infty }m_{N}^{\left( r\right) }=m^{\left( r\right) }.
\end{equation*}%
Then $F_{N}\rightarrow F$ vaguely.
\end{proposition}

See Theorem 4.5.5.on page 99 in \cite{chung01} for a proof. Note that Chung
uses words \textquotedblleft vague convergence\textquotedblright\ to denote
that kind of convergence which is more often called the weak convergence of
probability measures.

Since the semicircle distribution is bounded and therefore is determined by
its moments (see Corollary to Theorem II.12.7 in \cite{shiryaev95}),
therefore the assumption of Proposition \ref{weak_convergence} is satisfied,
and we only need to show that the moments of $S_{n}$ converge to the
corresponding moments of the semicircle distribution.

\textbf{Proof of Theorem \ref{CLT}:} Define $\eta _{i}$ as a sequence of
random variables that are freely independent among themselves and also
freely independent from all $\xi _{i}.$ Suppose also that $\eta _{i}$ have
semicircle distributions with $E\eta _{i}=0$ and $E\eta _{i}^{2}=\sigma
_{i}^{2}.$ We are going to accept the fact that the sum of free semicircle
random variables is semicircle, and therefore, the spectral distribution of $%
\left( \eta _{1}+...+\eta _{N}\right) /\left( s\sqrt{N}\right) $ converges
in distribution to the semicircle law $\mu _{SC}$ with zero expectation and
unit variance. Let us define $X_{i}=\xi _{i}/s_{N}$ and $Y_{i}=\eta
_{i}/s_{N}.$ We will proceed by proving that moments of $X_{1}+...+X_{N}$
converge to moments of $Y_{1}+...+Y_{N}$ and applying Proposition \ref%
{weak_convergence}. Let

\begin{equation*}
\Delta f=Ef\left( X_{1}+...+X_{N}\right) -Ef\left( Y_{1}+...+Y_{N}\right) ,
\end{equation*}%
where $f\left( x\right) =x^{m}.$ We want to show that this difference
approaches zero as $N$ grows.

By assumption, $EY_{i}=EX_{i}=0$ and $EY_{i}^{2}=EX_{i}^{2}=\sigma
_{i}^{2}/s_{N}^{2}.$

The first step is to write the difference $\Delta f$ as follows:%
\begin{eqnarray*}
\Delta f &=&\left[ Ef\left( X_{1}+...+X_{N-1}+X_{N}\right) -Ef\left(
X_{1}+...+X_{N-1}+Y_{N}\right) \right] \\
&&+\left[ Ef\left( X_{1}+...+X_{N-1}+Y_{N}\right) -Ef\left(
X_{1}+...+Y_{N-1}+Y_{N}\right) \right] \\
&&+\left[ Ef\left( X_{1}+Y_{2}+...+Y_{N-1}+Y_{N}\right) -Ef\left(
Y_{1}+Y_{2}+...+Y_{N-1}+Y_{N}\right) \right] .
\end{eqnarray*}%
We intend to estimate every difference in this sum. Let 
\begin{equation}
Z_{k}=X_{1}+...+X_{k-1}+Y_{k+1}+...+Y_{N}.  \label{Zk}
\end{equation}%
We are interested in 
\begin{equation*}
Ef\left( Z_{k}+X_{k}\right) -Ef\left( Z_{k}+Y_{k}\right) .
\end{equation*}

We are going to apply the Taylor expansion formula but first we define
directional derivatives. Let $f_{X_{k}}^{\prime }\left( Z_{k}\right) $ be
the \emph{derivative of }$f$\emph{\ at }$Z_{k}$\emph{\ in direction }$X_{k},$
defined as follows: 
\begin{equation*}
f_{X_{k}}^{\prime }(Z_{k})=:\lim_{t\downarrow 0}\frac{f\left(
Z_{k}+tX_{k}\right) -f(Z_{k})}{t}.
\end{equation*}%
The higher order directional derivatives can be defined recursively. For
example, 
\begin{equation}
f_{X_{k}}^{\prime \prime }\left( Z_{k}\right) =:\left( f_{X_{k}}^{\prime
}\right) _{X_{k}}^{\prime }\left( Z_{k}\right) =\lim_{t\downarrow 0}\frac{%
f_{X_{k}}^{\prime }\left( Z_{k}+tX_{k}\right) -f_{X_{k}}^{\prime }(Z_{k})}{t}%
.  \label{formula_second_derivative1}
\end{equation}%
For polynomials, this definition is equivalent to the following definition:%
\begin{equation}
f_{X_{k}}^{\prime \prime }(Z_{k})=2\lim_{t\downarrow 0}\frac{f\left(
Z_{k}+tX_{k}\right) -f(Z_{k})-tf_{X_{k}}^{\prime }(Z_{k})}{t^{2}}.
\label{formula_second_derivative2}
\end{equation}

\begin{example}
Operator directional derivatives of $f\left( x\right) =x^{4}$
\end{example}

Let us compute $f_{X}^{\prime }\left( Z\right) $ and $f_{X}^{\prime \prime
}\left( Z\right) $ for $f\left( x\right) =x^{4}.$ Using definitions we get 
\begin{equation*}
f_{X}^{\prime }\left( Z\right) =Z^{3}X+Z^{2}XZ+ZXZ^{2}+XZ^{3}
\end{equation*}%
and 
\begin{equation}
f_{X}^{\prime \prime }\left( Z\right) =2\left(
Z^{2}X^{2}+ZXZX+XZ^{2}X+ZX^{2}Z+XZXZ+X^{2}Z^{2}\right) ,
\label{formula_2nd_derivative_x4}
\end{equation}%
and the expression for $f_{X}^{\prime \prime }\left( Z\right) $ does not
depend on whether definition (\ref{formula_second_derivative1}) or (\ref%
{formula_second_derivative2}) was applied.

\bigskip

The derivatives of $f$ at $Z_{k}+\tau X_{k}$ in direction $X_{k}$ are
defined similarly, for example: 
\begin{eqnarray*}
&&f_{X_{k}}^{\prime \prime \prime }\left( Z_{k}+\tau X_{k}\right) \\
&=&6\lim_{t\downarrow 0}\frac{f\left( Z_{k}+\left( \tau +t\right)
X_{k}\right) -f(Z_{k}+\tau X_{k})-tf_{X_{k}}^{\prime }(Z_{k}+\tau X_{k})-%
\frac{1}{2}t^{2}f_{X_{k}}^{\prime \prime }(Z_{k}+\tau X_{k})}{t^{3}}.
\end{eqnarray*}

Next, let us write the Taylor formula for $f\left( Z_{k}+X_{k}\right) $:%
\begin{equation}
f\left( Z_{k}+X_{k}\right) =f(Z_{k})+f_{X_{k}}^{\prime }(Z_{k})+\frac{1}{2}%
f_{X_{k}}^{\prime \prime }\left( Z_{k}\right) +\frac{1}{2}\int_{0}^{1}\left(
1-\tau \right) ^{2}f_{X_{k}}^{\prime \prime \prime }\left( Z_{k}+\tau
X_{k}\right) d\tau .  \label{formula_Taylor_for_noncommut_polynomials}
\end{equation}%
Formula (\ref{formula_Taylor_for_noncommut_polynomials}) can be obtained by
integration by parts from the expression 
\begin{equation*}
f\left( Z_{k}+X_{k}\right) -f(Z_{k})=\int_{0}^{1}f_{X_{k}}^{\prime }\left(
Z_{k}+\tau X_{k}\right) d\tau .
\end{equation*}

For polynomials it is easy to write the explicit expressions for $%
f_{X_{k}}^{\left( r\right) }\left( Z_{k}\right) $ or $f_{X_{k}}^{\left(
r\right) }\left( Z_{k}+\tau X_{k}\right) $ although they can be quite
cumbersome for polynomials of high degree. Very schematically, for a
function $f\left( x\right) =x^{m},$ we can write 
\begin{equation}
f_{X_{k}}^{\prime }\left( Z_{k}\right)
=X_{k}Z_{k}^{m-1}+Z_{k}X_{k}Z_{k}^{m-2}+...+Z_{k}^{m-1}X_{k},
\label{first_derivative}
\end{equation}%
and 
\begin{equation}
f_{X_{k}}^{\prime \prime }\left( Z_{k}\right) =2\left(
X_{k}^{2}Z_{k}^{m-2}+X_{k}Z_{k}X_{k}Z_{k}^{m-3}+...+Z_{k}^{m-2}X_{k}^{2}%
\right) ,  \label{second_derivative}
\end{equation}%
Similar formulas hold for $f_{Y_{k}}^{\prime }\left( Z_{k}\right) $ and $%
f_{Y_{k}}^{\prime \prime }\left( Z_{k}\right) ,$ with the change that $Y_{k}$
should be used instead of $X_{k}.$

Using the assumptions that sequence $\left\{ X_{k}\right\} $ satisfies
Condition A and that variables $Y_{k}$ are free, we can conclude that $%
Ef_{Y_{k}}^{\prime }\left( Z_{k}\right) =Ef_{X_{k}}^{\prime }\left(
Z_{k}\right) =0$ and that $Ef_{Y_{k}}^{\prime \prime }\left( Z_{k}\right)
=Ef_{X_{k}}^{\prime \prime }\left( Z_{k}\right) .$ Indeed, consider, for
example, (\ref{second_derivative}). We can use expression (\ref{Zk}) for $%
Z_{k}$ and the free independence of $Y_{i}$ to expand (\ref%
{second_derivative}) as 
\begin{equation}
Ef_{X_{k}}^{\prime \prime }\left( Z_{k}\right) =\sum_{\alpha }c_{\alpha
}P_{\alpha }\left( E\left( X_{k}\overline{X_{1}}X_{k}\overline{X_{2}}\right)
,E\left( X_{k}\overline{X_{3}}X_{k}\overline{X_{4}}\right) ,...\right) ,
\label{Esecond_derivative}
\end{equation}%
where $\overline{X_{i}}$ denotes certain monomials in variables $%
X_{1},...,X_{k-1}$ (i.e., $\overline{X_{i}}=X_{i_{1}}...X_{i_{p}}$ with $%
i_{k}\in \left\{ 1,...,k-1\right\} $), and where $\alpha $ indexes certain
polynomials $P_{\alpha }.$ In other words, using the free independence of $%
Y_{i}$ and $X_{i}$ we expand the expectations of polynomial $%
f_{X_{k}}^{\prime \prime }\left( Z_{k}\right) $ as a sum over polynomials in
joint moments of variables $X_{j}$ and $Y_{i}$ where $j=1,...,k$ and $%
i=k+1,...,N.$ By freeness, we can reduce the resulting expression so that
the moments in the reduced expression are either joint moments of variables $%
X_{j}$ or joint moments of variables $Y_{i}$ but never involve both $X_{j}$
and $Y_{i}$. Moreover, we can explictly calculate the moments of $Y_{i}$
(i.e., expectations of the products of $Y_{i}$) because their are mutually
free. The resulting expansion is (\ref{Esecond_derivative}).

Let us try to make this process clearer by an example. Suppose that $f\left(
x\right) =x^{4},$ $N=4,$ $k=2$ and $Z_{k}=Z_{2}=X_{1}+Y_{3}+Y_{4}.$ We aim
to compute $Ef_{X_{2}}^{\prime \prime }\left( Z_{2}\right) .$ Using formula (%
\ref{formula_2nd_derivative_x4}), we write:%
\begin{eqnarray*}
Ef_{X_{2}}^{\prime \prime }\left( Z_{2}\right) &=&2E\left(
Z_{2}^{2}X_{2}^{2}+...\right) \\
&=&2E\left( \left( X_{1}+Y_{3}+Y_{4}\right) ^{2}X_{2}^{2}+...\right) \\
&=&2\{E\left( X_{1}^{2}X_{2}^{2}\right) +E\left( X_{1}Y_{3}X_{2}^{2}\right)
+E\left( X_{1}Y_{4}X_{2}^{2}\right) \\
&&+E\left( Y_{3}X_{1}X_{2}^{2}\right) +E\left( Y_{3}^{2}X_{2}^{2}\right)
+E\left( Y_{3}Y_{4}X_{2}^{2}\right) \\
&&+E\left( Y_{4}X_{1}X_{2}^{2}\right) +E\left( Y_{4}Y_{3}X_{2}^{2}\right)
+E\left( Y_{4}^{2}X_{2}^{2}\right) +...\}.
\end{eqnarray*}%
Then, using the freeness of $Y_{3}$ and $Y_{4}$ and the facts that $E\left(
Y_{i}\right) =0$ and $E\left( Y_{i}^{2}\right) =\sigma _{i}^{2},$ we
continue as follows:%
\begin{equation*}
Ef_{X_{2}}^{\prime \prime }\left( Z_{2}\right) =2\{E\left(
X_{1}^{2}X_{2}^{2}\right) +\sigma _{3}^{2}E\left( X_{2}^{2}\right) +\sigma
_{4}^{2}E\left( X_{2}^{2}\right) +...\},
\end{equation*}%
which is the expression we wanted to obtain.

It is important to note that the coefficients $c_{\alpha }$ do not depend on
variables $X_{j}$ but only on $Y_{j},$ $j>k,$ and on the locations which $%
Y_{j}$ take in the expansion of $f_{X_{k}}^{\prime \prime }\left(
Z_{k}\right) .$ Therefore, we can substitute $Y_{k}$ for $X_{k}$ and develop
a similar formula for $Ef_{Y_{k}}^{\prime \prime }\left( Z_{k}\right) $: 
\begin{equation}
Ef_{Y_{k}}^{\prime \prime }\left( Z_{k}\right) =\sum_{\alpha }c_{\alpha
}P_{\alpha }\left( E\left( Y_{k}\overline{X_{1}}Y_{k}\overline{X_{2}}\right)
,E\left( Y_{k}\overline{X_{3}}Y_{k}\overline{X_{4}}\right) ,...\right) .
\label{formula_expectation_2nd_derivative}
\end{equation}%
\qquad \qquad

In the example above, we will have 
\begin{equation*}
Ef_{Y_{2}}^{\prime \prime }\left( Z_{2}\right) =2\{E\left(
X_{1}^{2}Y_{2}^{2}\right) +\sigma _{3}^{2}E\left( Y_{2}^{2}\right) +\sigma
_{4}^{2}E\left( Y_{2}^{2}\right) +...\}.
\end{equation*}%
Formula (\ref{formula_expectation_2nd_derivative}) is exactly the same as
formula (\ref{Esecond_derivative})\ except that all $X_{k}$ are substituted
with $Y_{k}$. Finally, using Condition A we obtain that for every $i$: 
\begin{eqnarray*}
E\left( Y_{k}\overline{X_{i}}Y_{k}\overline{X_{i+1}}\right) &=&E\left(
Y_{k}^{2}\right) E\left( \overline{X_{i}}\right) E\left( \overline{X_{i+1}}%
\right) \\
&=&E\left( X_{k}^{2}\right) E\left( \overline{X_{i}}\right) E\left( 
\overline{X_{i+1}}\right) \\
&=&E\left( X_{k}\overline{X_{i}}X_{k}\overline{X_{i+1}}\right) ,
\end{eqnarray*}%
and therefore $Ef_{Y_{k}}^{\prime \prime }\left( Z_{k}\right)
=Ef_{X_{k}}^{\prime \prime }\left( Z_{k}\right) .$

Consequently, 
\begin{eqnarray*}
&&Ef\left( Z_{k}+X_{k}\right) -Ef\left( Z_{k}+Y_{k}\right) \\
&=&\frac{1}{2}\int_{0}^{1}\left( 1-\tau \right) ^{2}Ef_{X_{k}}^{\prime
\prime \prime }\left( Z_{k}+\tau X_{k}\right) d\tau -\frac{1}{2}%
\int_{0}^{1}\left( 1-\tau \right) ^{2}Ef_{Y_{k}}^{\prime \prime \prime
}\left( Z_{k}+\tau Y_{k}\right) d\tau .
\end{eqnarray*}

Next, note that if $f$ is a polynomial, then $f_{X_{k}}^{\prime \prime
\prime }\left( Z_{k}+\tau X_{k}\right) $ is the sum of a finite number of
terms which are products of $Z_{k}+\tau X_{k}$ and $X_{k}.$ The number of
terms in this expansion is bounded by $C_{1}$, which depends only on the
degree $m$ of the polynomial $f.$

A typical term in the expansion looks like 
\begin{equation*}
E\left( Z_{k}+\tau X_{k}\right) ^{m-7}X_{k}^{3}\left( Z_{k}+\tau
X_{k}\right) ^{3}X_{k}.
\end{equation*}%
In addition, if we expand the powers of $Z_{k}+\tau X_{k}$, we will get
another expansion that has the number of terms bounded by $C_{2},$ where $%
C_{2}$ depends only on $m.$ A typical element of this new expansion is 
\begin{equation*}
E\left( Z_{k}^{m-7}X_{k}^{3}Z_{k}^{2}X_{k}^{2}\right) .
\end{equation*}%
Every term in this expansion has a total degree of $X_{k}$ not less than $3,$
and, correspondingly, a total degree of $Z_{k}$ not more than $m-3.$ Our
task is to show that as $n\rightarrow \infty ,$ these terms approach $0.$

We will use the following lemma to estimate each of the summands in the
expansion of $f_{X_{k}}^{\prime \prime \prime }\left( Z_{k}+\tau
X_{k}\right) $.

\begin{lemma}
\label{lemma_Cauchy_Schwartz_generalized}Let $X$ and $Y$ be self-adjoint.
Then 
\begin{eqnarray*}
&&\left\vert E\left( X^{m_{1}}Y^{n_{1}}...X^{m_{r}}Y^{n_{r}}\right)
\right\vert \\
&\leq &\left[ E\left( X^{2^{r}m_{1}}\right) \right] ^{2^{-r}}\left[ E\left(
Y^{2^{r}n_{1}}\right) \right] ^{2^{-r}}...\left[ E\left(
X^{2^{r}m_{r}}\right) \right] ^{2^{-r}}\left[ E\left( Y^{2^{r}n_{r}}\right) %
\right] ^{2^{-r}}.
\end{eqnarray*}
\end{lemma}

\textbf{Proof:} For $r=1,$ this is the usual Cauchy-Schwartz inequality for
traces:%
\begin{equation*}
\left\vert E\left( X^{m_{1}}Y^{n_{1}}\right) \right\vert ^{2}\leq E\left(
X^{2m_{1}}\right) E\left( Y^{2n_{1}}\right) .
\end{equation*}%
See, for example, Proposition I.9.5 on p. 37 in \cite{takesaki79}.

Next, we proceed by induction. We have two slightly different cases to
consider. Assume first that $r$ is even, $\ r=2s.$ Then, by the
Cauchy-Schwartz inequality, we have: 
\begin{eqnarray*}
&&\left\vert E\left( X^{m_{1}}Y^{n_{1}}...X^{m_{r}}Y^{n_{r}}\right)
\right\vert ^{2} \\
&\leq &E\left(
X^{m_{1}}Y^{n_{1}}...X^{m_{s}}Y^{n_{s}}Y^{n_{s}}X^{m_{s}}...Y^{n_{1}}X^{m_{1}}\right) E\left( Y^{n_{r}}X^{m_{r}}...Y^{n_{s+1}}X^{m_{s+1}}X^{m_{s+1}}Y^{n_{s+1}}...X^{m_{r}}Y^{n_{r}}\right)
\\
&=&E\left(
X^{2m_{1}}Y^{n_{1}}...X^{m_{s}}Y^{2n_{s}}X^{m_{s}}...Y^{n_{1}}\right)
E\left(
Y^{2n_{r}}X^{m_{r}}...Y^{n_{s+1}}X^{2m_{s+1}}Y^{n_{s+1}}...X^{m_{r}}\right) .
\end{eqnarray*}%
Applying the inductive hypothesis, we obtain:%
\begin{eqnarray*}
&&\left\vert E\left( X^{m_{1}}Y^{n_{1}}...X^{m_{r}}Y^{n_{r}}\right)
\right\vert ^{2} \\
&\leq &\left[ E\left( X^{2^{r}m_{1}}\right) \right] ^{2^{-r+1}}\left[
E\left( Y^{2^{r}n_{s}}\right) \right] ^{2^{-r+1}}\left[ E\left(
Y^{2^{r-1}n_{1}}\right) \right] ^{2^{-r+2}}...\left[ E\left(
X^{2^{r-1}m_{s}}\right) \right] ^{2^{-r+2}} \\
&&\times \left[ E\left( X^{2^{r}m_{s+1}}\right) \right] ^{2^{-r+1}}\left[
E\left( Y^{2^{r}n_{r}}\right) \right] ^{2^{-r+1}}\left[ E\left(
Y^{2^{r-1}n_{s+1}}\right) \right] ^{2^{-r+2}}...\left[ E\left(
X^{2^{r-1}m_{r}}\right) \right] ^{2^{-r+2}}.
\end{eqnarray*}%
We recall that by the Lyapunov inequality, $\left[ E\left(
Y^{2^{r-1}n_{1}}\right) \right] ^{2^{-r+2}}\leq \left[ E\left(
Y^{2^{r}n_{1}}\right) \right] ^{2^{-r+1}}$ and we get the desired inequality:%
\begin{eqnarray*}
&&\left\vert E\left( X^{m_{1}}Y^{n_{1}}...X^{m_{r}}Y^{n_{r}}\right)
\right\vert \\
&\leq &\left[ E\left( X^{2^{r}m_{1}}\right) \right] ^{2^{-r}}\left[ E\left(
Y^{2^{r}n_{1}}\right) \right] ^{2^{-r}}...\left[ E\left(
X^{2^{r}m_{r}}\right) \right] ^{2^{-r}}\left[ E\left( Y^{2^{r}n_{r}}\right) %
\right] ^{2^{-r}}.
\end{eqnarray*}

Now let $r$ be odd, $r=2s+1.$ Then 
\begin{eqnarray*}
&&\left\vert E\left( X^{m_{1}}Y^{n_{1}}...X^{m_{r}}Y^{n_{r}}\right)
\right\vert ^{2} \\
&\leq &E\left(
X^{m_{1}}Y^{n_{1}}...Y^{n_{s}}X^{m_{s+1}}X^{m_{s+1}}Y^{n_{s}}...Y^{n_{1}}X^{m_{1}}\right) E\left( Y^{n_{r}}X^{m_{r}}...X^{m_{s+2}}Y^{n_{s+1}}Y^{n_{s+1}}X^{m_{s+2}}...X^{m_{r}}Y^{n_{r}}\right)
\\
&=&E\left(
X^{2m_{1}}Y^{n_{1}}...Y^{n_{s}}X^{2m_{s+1}}Y^{n_{s}}...Y^{n_{1}}\right)
E\left(
Y^{2n_{r}}X^{m_{r}}...X^{m_{s+2}}Y^{2n_{s+1}}X^{m_{s+1}}...X^{m_{r}}\right) .
\end{eqnarray*}%
After that we can use the inductive hypothesis and the Lyapunov inequality
and obtain that 
\begin{eqnarray*}
&&\left\vert E\left( X^{m_{1}}Y^{n_{1}}...X^{m_{r}}Y^{n_{r}}\right)
\right\vert \\
&\leq &\left[ E\left( X^{2^{r}m_{1}}\right) \right] ^{2^{-r}}\left[ E\left(
Y^{2^{r}n_{1}}\right) \right] ^{2^{-r}}...\left[ E\left(
X^{2^{r}m_{r}}\right) \right] ^{2^{-r}}\left[ E\left( Y^{2^{r}n_{r}}\right) %
\right] ^{2^{-r}}.
\end{eqnarray*}%
QED.

We apply Lemma \ref{lemma_Cauchy_Schwartz_generalized} to estimate each of
the summands in the expansion of $f_{X_{k}}^{\prime \prime \prime }\left(
Z_{k}+\tau X_{k}\right) $. Consider a summand $E\left(
Z_{k}^{m_{1}}X_{k}^{n_{1}}...Z_{k}^{m_{r}}X_{k}^{n_{r}}\right) .$ Then by
Lemma \ref{lemma_Cauchy_Schwartz_generalized}, we have 
\begin{eqnarray}
&&\left\vert E\left(
Z_{k}^{m_{1}}X_{k}^{n_{1}}...Z_{k}^{m_{r}}X_{k}^{n_{r}}\right) \right\vert
\label{estimate_third_derivative} \\
&\leq &\left[ E\left( Z_{k}^{2^{r}m_{1}}\right) \right] ^{2^{-r}}\left[
E\left( X_{k}^{2^{r}n_{1}}\right) \right] ^{2^{-r}}...\left[ E\left(
Z_{k}^{2^{r}m_{r}}\right) \right] ^{2^{-r}}\left[ E\left(
X_{k}^{2^{r}n_{r}}\right) \right] ^{2^{-r}}.  \notag
\end{eqnarray}

Next step is to estimate the absolute moments of the variable $Z_{k}.$

\begin{lemma}
\label{moment_bound0}Let $Z=\left( v_{1}+...+v_{N}\right) /N^{1/2},$ where $%
v_{i}$ are self-adjoint and satisfy condition A(1) and let $E\left\vert
v_{i}\right\vert ^{k}\leq \mu _{k}$ for every $i.$ Then, for every integer $%
r\geq 0$ 
\begin{equation*}
E\left( \left\vert Z\right\vert ^{r}\right) =O\left( 1\right) \text{ as }%
N\rightarrow \infty .
\end{equation*}
\end{lemma}

\textbf{Proof:} We will first treat the case of even $r.$ In this case, $%
E\left( \left\vert Z\right\vert ^{r}\right) =E\left( Z^{r}\right) .$
Consider the expansion of $\left( v_{1}+...+v_{N}\right) ^{r}.$ Let us refer
to the indices $1,$ $...,$ $N$ as colors of the corresponding $v.$ If a term
in the expansion includes more than $r/2$ distinct colors, then one of the
colors must be used by this term only once. Therefore, by the first part of
condition A the expectation of such a term is 0.

Let us estimate a number of terms in the expansion that include no more than 
$r/2$ distinct colors. Consider a fixed combination of $\leq r/2$ colors.
The number of terms that use colors only from this combination is $\leq
\left( r/2\right) ^{r}.$ Indeed, consider the product \newline
$\left( v_{1}+...+v_{N}\right) \left( v_{1}+...+v_{N}\right) ...\left(
v_{1}+...+v_{N}\right) $ with $r$ product terms. We can choose an element
from the first product term in $r/2$ possible ways, an element from the
second product term in $r/2$ possible ways, etc. Therefore, the number of
all possible choices is $\left( r/2\right) ^{r}.$ On the other hand, the
number of possible different combinations of $k\leq r/2$ colors is 
\begin{equation*}
\frac{N!}{\left( N-k\right) !k!}\leq N^{r/2}.
\end{equation*}%
Therefore, the total number of terms that use no more than $r/2$ colors is
bounded from above by 
\begin{equation*}
\left( r/2\right) ^{r}N^{r/2}.
\end{equation*}

Now let us estimate the expectation of an individual term in the expansion.
In other words we want to estimate $E\left(
v_{i_{1}}^{k_{1}}...v_{i_{s}}^{k_{s}}\right) ,$ where $k_{t}\geq 1,$ $%
k_{1}+...+k_{s}=r,$ and $i_{t}\neq i_{t+1}.$ First, note that 
\begin{equation*}
\left\vert E\left( v_{i_{1}}^{k_{1}}...v_{i_{s}}^{k_{s}}\right) \right\vert
\leq E\left( \left\vert v_{i_{1}}^{k_{1}}...v_{i_{s}}^{k_{s}}\right\vert
\right) .
\end{equation*}%
Indeed, using the Cauchy-Schwartz inequality, for any operator $X$ we can
write 
\begin{eqnarray*}
\left\vert E\left( X\right) \right\vert ^{2} &=&\left\vert E\left(
U\left\vert X\right\vert ^{1/2}\left\vert X\right\vert ^{1/2}\right)
\right\vert ^{2}\leq E\left( \left\vert X\right\vert ^{1/2}U^{\ast
}U\left\vert X\right\vert ^{1/2}\right) E\left( \left\vert X\right\vert
^{1/2}\left\vert X\right\vert ^{1/2}\right) \\
&=&E\left( \left\vert X\right\vert P\right) E\left( \left\vert X\right\vert
\right) ,
\end{eqnarray*}%
where $U$ is a partial isometry and $P=U^{\ast }U$ is a projection. Note
that from the positivity of the expectation functional it follows that $%
E\left( \left\vert X\right\vert P\right) \leq E\left( \left\vert
X\right\vert \right) .$ Therefore, we can conclude that $\left\vert E\left(
X\right) \right\vert \leq E\left( \left\vert X\right\vert \right) .$

Next, we use the H\"{o}lder inequality for traces of non-commutative
operators (see \cite{fack82}, especially Corollary 4.4(iii) on page 324, for
the case of the trace in a von Neumann algebra and Section III.7.2 in \cite%
{gohberg_krein69} for the case of compact operators and the usual operator
trace). Note that 
\begin{equation*}
\underset{s\text{-times}}{\underbrace{\frac{1}{s}+...+\frac{1}{s}}}=1,
\end{equation*}%
therefore, the H\"{o}lder inequality gives%
\begin{equation*}
E\left( \left\vert v_{i_{1}}^{k_{1}}...v_{i_{s}}^{k_{s}}\right\vert \right)
\leq \left[ E\left( \left\vert v_{i_{1}}\right\vert ^{k_{1}s}\right)
...E\left( \left\vert v_{i_{s}}\right\vert ^{k_{s}s}\right) \right] ^{1/s}.
\end{equation*}

Using this result and the uniform boundedness of the moments (from
assumption of the lemma), we get:%
\begin{equation*}
\log \left\vert E\left( v_{i_{1}}^{k_{1}}...v_{i_{s}}^{k_{s}}\right)
\right\vert \leq \frac{1}{s}\sum_{i=1}^{s}\log \mu _{k_{i}s}.
\end{equation*}%
Without loss of generality we can assume that bounds $\mu _{k}$ are
increasing in $k.$ Using the facts that $s\leq r$ and $k_{i}\leq r$, we
obtain the bound:%
\begin{equation*}
\log \left\vert E\left( v_{i_{1}}^{k_{1}}...v_{i_{s}}^{k_{s}}\right)
\right\vert \leq \log \mu _{r^{2}},
\end{equation*}%
or 
\begin{equation*}
\left\vert E\left( v_{i_{1}}^{k_{1}}...v_{i_{s}}^{k_{s}}\right) \right\vert
\leq \mu _{r^{2}}.
\end{equation*}

Therefore,%
\begin{equation*}
E\left( v_{1}+...+v_{N}\right) ^{r}\leq \left( r/2\right) ^{r}\mu
_{r^{2}}N^{r/2},
\end{equation*}%
and 
\begin{equation}
E\left( Z^{r}\right) \leq \left( r/2\right) ^{r}\mu _{r^{2}}.
\label{inequality_bound_even_moments}
\end{equation}

Now consider the case of odd $r.$ In this case, we use the Lyapunov
inequality to write:%
\begin{eqnarray}
E\left\vert Z\right\vert ^{r} &\leq &\left( E\left\vert Z\right\vert
^{r+1}\right) ^{\frac{r}{r+1}}  \label{inequality_bound_odd_moments} \\
&\leq &\left( \left( \frac{r+1}{2}\right) ^{r+1}\mu _{\left( r+1\right)
^{2}}\right) ^{\frac{r}{r+1}}  \notag \\
&=&\left( \frac{r+1}{2}\right) ^{r}\left( \mu _{\left( r+1\right)
^{2}}\right) ^{\frac{r}{r+1}}.  \notag
\end{eqnarray}%
The important point is that the bounds in (\ref%
{inequality_bound_even_moments}) and (\ref{inequality_bound_odd_moments}) do
not depend on $N.$ QED.

By definition $Z_{k}=\left( \xi _{1}+...+\xi _{k-1}+\eta _{k+1}+...+\eta
_{N}\right) /s_{N}$ and by assumption $\xi _{i}$ and $\eta _{i}$ are
uniformly bounded, and $s_{N}\sim \sqrt{N}$. Moreover, $\xi _{1},...,\xi
_{k-1}$ satisfy Condition A by assumption, and $\eta _{k+1},...,\eta _{N}$
are free from each other and from $\xi _{1},...,\xi _{k-1}.$ Therefore, by
Proposition \ref{proposition_A1}, $\xi _{1},...,\xi _{k-1},\eta
_{k+1},...,\eta _{N}$ satisfy condition A(1). Consequently, we can apply
Lemma \ref{moment_bound0} to $Z_{k}$ and conclude that $E\left\vert
Z_{k}\right\vert ^{r}$ is bounded by a constant that depends only on $r$ but
does not depend on $N.$

Using this fact, we can continue the estimate in (\ref%
{estimate_third_derivative}) and write: 
\begin{eqnarray}
&&\left\vert E\left(
Z_{k}^{m_{1}}X_{k}^{n_{1}}...Z_{k}^{m_{r}}X_{k}^{n_{r}}\right) \right\vert \\
&\leq &C_{4}\left[ E\left( X_{k}^{2^{r}n_{1}}\right) \right] ^{2^{-r}}...%
\left[ E\left( X_{k}^{2^{r}n_{r}}\right) \right] ^{2^{-r}},  \notag
\end{eqnarray}
where the constant $C_{4}$ depends only on $m.$

Next we note that 
\begin{equation*}
\left[ E\left( X_{k}^{2^{r}n_{1}}\right) \right] ^{2^{-r}}\leq C\left( \frac{%
\mu _{2^{r}n_{1}}}{N^{2^{r-1}n_{1}}}\right) ^{2^{-r}}=C\frac{\left( \mu
_{2^{r}n_{1}}\right) ^{2^{-r}}}{N^{n_{1}/2}}.
\end{equation*}%
Next note that $n_{1}+...+n_{r}\geq 3;$ therefore we can write 
\begin{equation*}
\left[ E\left( X_{k}^{2^{r}n_{1}}\right) \right] ^{2^{-r}}...\left[ E\left(
X_{k}^{2^{r}n_{r}}\right) \right] ^{2^{-r}}\leq C^{\prime }N^{-3/2}.
\end{equation*}

In sum$,$ we obtain the following Lemma:

\begin{lemma}
\begin{equation*}
\left\vert Ef_{X_{k}}^{\prime \prime \prime }\left( Z_{k}+\tau X_{k}\right)
\right\vert \leq C_{5}N^{-3/2},
\end{equation*}%
where $C_{5}$ depends only on the degree of polynomial $f$ and the sequence
of constants $\mu _{k}$.
\end{lemma}

A similar result holds for $\left\vert Ef_{X_{k}}^{\prime \prime \prime
}\left( Z_{k}+\tau Y_{k}\right) \right\vert $ and we can conclude that 
\begin{equation*}
\left\vert Ef\left( Z_{k}+X_{k}\right) -Ef\left( Z_{k}+Y_{k}\right)
\right\vert \leq C_{6}N^{-3/2}.
\end{equation*}%
After we add these inequalities over all $k=1,...,N$ we get 
\begin{equation*}
\left\vert Ef\left( X_{1}+...+X_{N}\right) -Ef\left( Y_{1}+...+Y_{N}\right)
\right\vert \leq C_{7}N^{-1/2}.
\end{equation*}%
Clearly this estimate approaches $0$ as $N$ grows. Applying Proposition \ref%
{weak_convergence}, we conclude that the measure of $X_{1}+...+X_{N}$
converges to the measure of $Y_{1}+...+Y_{N}$ in distribution. This finishes
the proof of the main theorem.

\section{Concluding Remarks}

The key points of this proof are as follows: 1) We can substitute each
random variable $X_{i}$ in the sum $S_{N}$ with a free random variable $%
Y_{i} $ so that the first and the second derivatives of any polynomial with $%
S_{N}$ in the argument remain unchanged. The possibility of this
substitution depends on Condition A being satisfied by $X_{i}.$ 2)We can
estimate a change in the third derivative as we substitute $Y_{i}$ for $%
X_{i} $ by using the first part of Condition A and several matrix
inequalities, valid for any collection of operators. Here Condition A is
used only in the proof that the $k$-th moment of $\left( \xi _{1}+...+\xi
_{N}\right) /N^{1/2} $ is bounded as $N\rightarrow \infty $.

It is interesting to speculate whether the ideas in this proof can be
generalized to the case of the multivariate CLT.

\end{document}